\documentclass{article}[12pt]
\usepackage{amsfonts,amsmath}

\usepackage{color}

\newtheorem{theorem}{Theorem}
\newtheorem{proposition}{Proposition}

\def\C{{\mathbb C}}

\begin{document}

\title{Moutard type transformation for matrix generalized analytic functions
and gauge transformations
\thanks{The work was done during the visit of the second author (I.A.T.) to Centre de Math\'ematiques Appliqu\'ees of 
\'Ecole Polytechnique in July 2016 and was supported by
RSF (grant 14-11-00441).}}

\author{R.G. Novikov \thanks{CNRS (UMR 7641), Centre de Math\'ematiques
Appliqu\'ees, \'Ecole Polytechnique, 91128 Palaiseau, France, and
IEPT RAS, 117997 Moscow, Russia; e-mail:
novikov@cmap.polytechnique.fr} \and I.A. Taimanov
\thanks{Sobolev Institute of Mathematics, 630090 Novosibirsk,
Russia, and Novosibirsk State University, 630090 Novosibirsk, Russia; e-mail: taimanov@math.nsc.ru}}

\date{}

\maketitle

Considerable progress in the theory of Darboux-Moutard type transformations for two-dimensional linear differential systems with 
applications to geometry, spectral theory, and soliton equations has been achieved recently, see, e.g., \cite{L,TT,T,GN}. 
In the present note we derive such a transformation for the matrix generalized function system
\begin{equation}
    \label{1}
\partial_{\bar{z}}\Psi + A\Psi + B \bar{\Psi} = 0,
\end{equation}
where  $\partial_{\bar{z}} = \frac{\partial}{\partial \bar{z}}$, the coefficients $A$ and $B$ and solutions $\Psi$ are 
$(N \times N)$-matrix functions on $D$, with $D$ an open simply connected domain in $\C$.
In particular, this generalizes the transformation for $N=1$ found in \cite{GN} with $A=0$.
In addition, we show that the Moutard type transformation for system
(\ref{1}) with $B=0$ is equivalent to a gauge transformation for the connection $\nabla_{\bar{z}} = \partial_{\bar{z}} + A$.
In turn, our studies show that the Moutard type transformation for system
(\ref{1}) with $A=0$ 
can be treated as a proper analog of the forementioned gauge transformation.  

As for $N=1$, system (\ref{1}) is reduced to the system 
\begin{equation}
\label{2}
\partial_{\bar{z}}\Psi + B \bar{\Psi} = 0,    
\end{equation}
i.e. to system (\ref{1}) with $A=0$,
by the gauge transformation
$$
\Psi \to \tilde{\Psi}=g^{-1}\Psi, \ \ B \to \tilde{B}=g^{-1}B\bar{g}, \ \ \partial_{\bar{z}} g + Ag=0, \  \ \det g \neq 0.
$$
We say that the system
\begin{equation}
\label{3}
\partial_z \Psi^+ - \bar{\Psi}^+ B = 0   
\end{equation}
is conjugate to system (\ref{2}) (see \cite{V} for a similar definition for $N=1$).

We have the following result.

\begin{theorem}
\label{t1}
Systems (\ref{2}) and (\ref{3}) are covariant, i.e. mapped into the systems of the same type, with respect to the Moutard type transformation
\begin{equation}
\label{4} 
\begin{split}
\Psi \to \tilde{\Psi} = \Psi - F \,\omega_{F,F^+}^{-1}\,\omega_{\Psi,F^+}, \\
\Psi^+ \to \tilde{\Psi}^+ = \Psi^+ - \omega_{F,\Psi^+}\,\omega_{F,F^+}^{-1}\, F^+, \\ 
B \to \tilde{B}=B + F \,\omega_{F,F^+}^{-1}\, F^+,
\end{split}
\end{equation}
where $F$ and $F^+$ are arbitrary fixed solutions of (\ref{2}) and (\ref{3}), respectively, 
\begin{equation}
    \label{5}
\partial_{\bar{z}} \omega_{\Phi,\Phi^+} = \Phi^+ \bar{\Phi}, \ \ \
\mathrm{Re}\, \omega_{\Phi,\Phi^+} = 0, 
\end{equation}
for $\Phi$ and $\Phi^+$ meeting equations (\ref{2}) and (\ref{3}), and
$\det\,\omega_{F,F^+} \neq 0$.
\end{theorem}

For finding $\omega_{\Phi,\Phi^+}$ satisfying (\ref{5}) we use also that $\partial_z \omega_{\Phi,\Phi^+} = 
-\bar{\Phi}^+ \Phi$. In addition, our definition of $\omega_{\Phi,\Phi^+}$ is self-consistent up to 
a pure imaginary matrix integration constant in view of the identity $\partial_z \Phi^+ \bar{\Phi} = - \partial_{\bar{z}} 
\bar{\Phi}^+ \Phi$. The latter equality follows from systems (\ref{2}) and (\ref{3}) for $\Phi$ and $\Phi^+$, respectively.  
We recall that the domain $D$ is simply connected. 

Given $\omega_{F,F^+}, \omega_{\Psi,F^+}$, and $\omega_{F,\Psi^+}$, 
Theorem \ref{t1} is proved by straightforward computations.

In addition, for the system
\begin{equation}
    \label{6}
\partial_{\bar{z}}\Psi + A\Psi = 0,
\end{equation}
i.e., for system (\ref{1}) with $B=0$,
the following result also holds.

\begin{proposition}
\label{p1}
System (\ref{6})  is covariant under the following Moutard type transformation
\begin{equation}
    \label{7}
\Psi \to \tilde{\Psi} = \Psi - F \,\hat{\omega}_{F,F^+}^{-1}\,\hat{\omega}_{\Psi,F^+}, \ \
A \to \tilde{A}= A + F \,\hat{\omega}_{F,F^+}^{-1}\, F^+,
\end{equation}
 where $F$ is an arbitrary fixed solution of (\ref{6}), $F^+$ is an arbitrary fixed matrix function, 
\begin{equation}
  \label{8}
\partial_{\bar{z}} \hat{\omega}_{\Phi,F^+} = F^+\Phi
\end{equation}
for any matrix function $\Phi$, and
$\det\, \hat{\omega}_{F,F^+} \neq 0$.
\end{proposition}

Equations (\ref{7}) and (\ref{8}) are analogs of equations (\ref{4}) and (\ref{5}).
However, in difference with (\ref{5}), we do not require that the matrix functions 
$\hat{\omega}_{F,F^+}$ would be pure imaginary. Equation (\ref{8}) is solvable for $\hat{\omega}_{F,F^+}$ 
and Proposition 1 is proved by straightforward computations. 

{\sc Remark.}
Let $A, \tilde{A}, \Psi, F, F^+$, and $\hat{\omega}_{\Phi,F^+}$ be the same as in Proposition 1. Let 
$$ 
g = 1 - F \hat{\omega}_{F,F^+}^{-1}\Lambda, \ \ \Lambda_{\bar{z}} = \Lambda A + F^+.
$$
Then
$$
\partial_{\bar{z}} (g\Psi) + \tilde{A} (g\Psi) =0.    
$$
It is proved by straightforward computations and it shows that for invertible $g$ 
the transformation $A \to \tilde{A}$ reduces to a gauge transform.

\end{document}